# ESSENTIALLY CRITICALLY COMPRESSIBLE MODULES

Abhay K. Singh[1]

[1] Department of Applied Mathematics, Indian School of Mines, Dhanbad-826004,India

**ABSTRACT**

The main purpose of this paper is to introduce the concept of essentially critically compressible modules. We call an R-module M essentially critically compressible module if it is essentially compressible and additionally it cannot be imbedded in any of its factor module $M/N$, where N is essential submodule of M. **Keywords: template, Word**

## 1 INTRODUCTION

Throughout all rings have non-zero identity elements and all modules are unital right R-modules. J. M. Zelmanowitz [1976] defined, an R-module M to be compressible if it can be imbedded in any of its non-zero submodules. A compressible module is called critically compressible if it cannot be imbedded in any of its proper factor modules. Zelmanowitz [1977] proved the following result about critically compressible modules. Zelmanowitz [1981] introduced weakly primitive rings and proved that a ring is weakly primitive if and only if it is a weakly dense subring of a full linear ring. This is a generalization of the Jacobson density theorem. J.M. Zelmanowitz [1993] gave the notion of weakly compressible modules. According to Zelmanowitz, an R-module M is called weakly compressible if for a non-zero submodule N of M, there exists a non-zero homomorphism f from M to N such that $f^2$ is non-zero. P.F. Smith [2005] introduced the concept of a slightly compressible module, which is a generalization of compressible and weakly compressible modules. According to Smith, the module M is called slightly compressible if, for a non-zero submodule N of M, there exists a non-zero homomorphism from M to N. With the help of an example, Smith has shown that for every ring a projective module is not necessarily slightly compressible, and proved that over a V-ring R, every projective R-module is slightly compressible. P.F. Smith [2006] introduced the concept of an essentially compressible module, which is another generalization of compressible modules. According to Smith, an R-module M is called essentially compressible module if, for a non-zero essential submodule N of M, there exists a non-zero monomorphism from M to N. Smith proved that a non-singular essentially compressible module M is isomorphic to a submodule of a free module, and the converse holds if R is a semiprime right Goldie ring. In case, R is a right FBN ring, M is essentially compressible if and only if M is subisomorphic to a critically compressible module. Abhay K. Singh [2007] discussed the characterization of a critically compressible module is terms of torsion free and uniform modules. V. S. Rodrigues and A. A. Sant'Ana[2009] considerd a problem due to Zelmanowitz. Specifically, they study under what conditions a uniform compressible modules whose non-zero endomorphisms are monomorphisms are critically compressible and also provided positive answers. Abhay K. Singh [2011] introduced and investigated the concept of essentially slightly compressible modules and rings, which are natural generalizations of essentially compressible modules and rings. He has also provided an example of an essentially slightly compressible module which is not essentially compressible. Cesim Celik [2012] introduced and investigated the concept of completely slightly compressible modules as a generalization of compressible modules.

**2. Essentially critically compressible Modules:**

**Definition2.1.** Let $M$ be a non-zero right $R-$ module. Then:
(1) $M$ is called essentially compressible if it can be imbedded in each of its non-zero essential submodules.
(2) $M$ is called essentially critically compressible if it is essentially compressible and, in addition, it cannot be imbedded in any of its factor modules $M/N$, where $N$ is essential submodule of $M$.

**Definition2.2.** An essentially partial endomorphism of a module $M$ is a homomorphism from an essential submodule of $M$ in to $M$.
Clearly for a uniform module $M$, every essentially partial endomorphism is a partial endomorphism on $M$.

**Example 2.1.** Clearly every critically compressible module is essentially critically compressible. However, every semisimple module over commutative ring is essentially critically compressible but need not be critically compressible.

We note that, every non-zero essential submodule of essentially compressible or essentially critically compressible module is also essentially compressible or essentially critically compressible module respectively.

There is useful theorem concerning critically compressible module, which was prompted in [1981].

**Proposition1.1.** [Zelmanowitz (1981), Prop.1.1] Let $M$ be a compressible module. Then the following conditions are equivalent:

(1) M is critically compressible;
(2) Every non-zero partial endomorphism on M is monomorphism.

Now above result can be extended for essentially compressible module $M$.

**Proposition1.2.** The following conditions are equivalent for a essentially compressible module M:

(1) M is essentially critically compressible;
(2) Every non-zero essentially partial endomorphism on M is monomorphism.

Proof. (1) → (2)

Let $N$ be a non-zero essential submodule of $M$ and let $f : N → M$ be a non-zero essential partial endomorphism on $M$. Then we have an isomorphism $\overline{f: N/_{Kerf}} \to f(N)$, and $f(N)$ is a non-zero essential submodule of $M$ [Anderson and fuller]. Since $M$ is essentially compressible, there is a monomorphism $g: M \to f(N)$ and the composition

$$M \xrightarrow{g} f(N) \xrightarrow{f^{-1}} N/_{Kerf} \leq M/_{Kerf}$$

is a monomorphism. By hypothesis $Kerf = 0$.

(2) → (1)

Assume $M$ is not essentially critically compressible module. Then there is a non-zero proper essential submodule $N$ of $M$ is imbedded in $M/_N$, say $h: M \to M/_N$ is the monomorphism. Let $T$ be the submodule of M such that $M \leq T \leq N$. Then we have, the composition

$$T \xrightarrow{\pi} T/_N \xrightarrow{h^{-1}} M$$

is non-zero essential partial endomorphism and so monomorphism, where $\pi$ is the canonical projection. But $N \neq 0$ is contained in the kernel of the composition.

This is a contradiction.

An $R$ − module $M$ is called essentially retractable [2007] if $Hom_R(M, N) \neq 0$ for every essential submodule $N$ of $M$. It is clear that every essentially compressible module is essentially retractable but converse is not true. Indeed, $Z_4$ over $Z$ is essentially retractable but not essentially compressible.

**Proposition 1.3.** Suppose that $M$ is an essentially retractable $R$-module. If every non-zero $f \epsilon End(M)$ is a monomorphism, then $M$ is essentially compressible.

Proof. Let N be a non-zero essential submodule of $M$ and $f : M \to N$ a non-zero homomorphism. Then $fi$ is a monomorphism and obviously $g$ is a monomorphism, where $i$ is an inclusion map. Converse is obvious.

Then Proposition 1.2 can be extended to the setting of essentially retractable modules.

**Proposition 1.4.** Let $M$ be an essentially retractable $R$-module. The following statements are equivalent:

(i) M is essentially critically compressible;

(ii) Every non-zero essential partial endomorphism on $M$ is monomorphism.

Now the following results give another formulation to the Zelmanowitz's question in terms of essential compressible modules.

**Theorem 1.5.** Let $M$ be an $R$-module. The following conditions are equivalent:

(i) $M$ is compressible and every non-zero endomorphism of $M$ is monomorphism;

(ii) $M$ is essentially compressible and $End(M)$ is domain;

(iii) $M$ is essentially retractable and every non-zero endomorphism of $M$ is monomorphism;

(iv) $M$ is essentially retractable and $End(M)$ is domain.

Proof. It is clear that (i) ⇒(ii), (ii) ⇒ (iv) and (iv) ⇒ (iii) and (iii) ⇒ (i) follows from Proposition 1.3.

**Proposition 1.6.** Let $M$ be a essentially retractable uniform module such that $End(M)$ is domain. The following conditions are equivalent:

(i) $M$ is critically compressible;
(ii) $M$ is essentially critically compressible;
(iii) $M$ is polyform.

Proof. (i) ⇔ (ii) For a uniform module critically compressible and essentially critically compressible modules are equivalent.

(iii) ⇒ (i). Suppose that $M$ is polyform. Since a module is polyform and uniform if and only if it is monoform [ ], by using that $M$ essentially retractable, we have that $M$ is essentially critically compressible.

(ii)⇒ (iii). If $M$ is essentially critically compressible, then due to uniformity of $M$, it is critically compressible. Then by Proposition 1.4 it is monoform, and hence polyform.

**Definition2.3.** A module $M$ is said to be essentially fully retractable if for every non-zero essential submodule $N$ of $M$ and every non-zero element $f \epsilon Hom_R(N, M)$ we have $Hom_R(M, N)f \neq 0$.

Clearly, if $M$ is essentially fully retractable then $M$ is retractable. A non-zero $R$-module $M$ is called self-similar if every non-zero submodule of $M$ is isomorphic to $M$. It is clear that self-similar modules are essentially fully retractable, but the converse is not true. For example $Z_4$ is a essentially fully retractable $Z$-module which is not self-similar.

**Proposition 1.7.** Let $M$ be essentially retractable uniform module such that $End(M)$ is domain. Then following conditions are equivalent:

(i) M is critically compressible;
(ii) M is essentially critically compressible;
(iii) M is polyform;
(iv) M is essentially fully retractable.

Proof. (i)⇔(ii) follows directly due to uniformity of $M$. (i)⇔(iii) follows directly from Proposition 1.6. Now we prove (iv)⇒(iii). Suppose that $M$ is not polyform. Then there exist a non-zero essential submodule $N$ of $M$ and a non-zero homomorphism $f: N \to M$ such that $Ker(f)$ is essential in $N$. But we have that $Hom_R(M,N)f \neq 0$. So there exist a non-zero map $g: M \to N$ such that $gf \neq 0$ and it follows that $gf$ is a monomorphism (see Theorem 1.5). Thus $g$ is monomorphism. Now $0 \neq Ker(gf) = g^{-1}(Ker(f)) \cong Ker(f) \cap Img(g)$, because $g$ is monomorphism. Since $Ker(f)$ is essential in $N$, we have necessarily $Im(g) = 0$ which is a contradiction.

(iii)⇒(iv).Since $M$ is polyform and uniform, then it is monoform. Therefore if $N$ is non-zero essential submodule of $M$ and $f$ is non-zero homomorphism from $N$ to $M$, then it is monomorphism. Since $M$ is essentially retractable, $Hom_R(M,N)$ is non-zero and we have that $Hom_R(M,N)f \neq 0$.

## ACKNOWLEDGEMENTS

This work is partially supported by ISM Dhanbad under the project FRS(15)/2010-2011/AM